\documentclass[a4paper,10pt]{article}
\usepackage{stmaryrd}
\usepackage{amsfonts}
\usepackage{bbm}
\usepackage{amscd}
\usepackage{mathrsfs}
\usepackage{latexsym,amssymb,amsmath,amscd,amscd,amsthm,amsxtra,xypic}
\usepackage[dvips]{graphicx}

\newtheorem{thm}{Theorem}[section]
\newtheorem{lem}[thm]{Lemma}
\newtheorem{cor}[thm]{Corollary}
\newtheorem{pro}[thm]{Proposition}
\newtheorem{ex}[thm]{Example}
\newtheorem{rmk}[thm]{Remark}
\newtheorem{defi}[thm]{Definition}

\newcommand{\lon }{\,\rightarrow\,}
\newcommand{\be }{\begin{eqnarray*}}
\newcommand{\ee }{\end{eqnarray*}}

\newcommand{\defbe}{\triangleq}
\newcommand{\pf}{\noindent{\bf Proof.}\ }

\newcommand{\inverse}{^{-1}}

\newcommand{\Real}{\mathbb R}

\newcommand{\CWM}{C^{\infty}(M)}

\newcommand{\set}[1]{\left\{#1\right\}}

\newcommand{\pairing}[1]{\left\langle #1\right\rangle}

\newcommand{\frkd}{\mathfrak d}
\newcommand{\frke}{\mathfrak e}

\newcommand{\frkp}{\mathfrak p}

\newcommand{\frkC}{\mathfrak C}
\newcommand{\frkD}{\mathfrak D}

\newcommand{\frkX}{\mathfrak X}

\def\gpd{\,\lower1pt\hbox{$\longrightarrow$}\hskip-.24in\raise2pt
         \hbox{$\longrightarrow$}\,}

\def\qed{\hfill ~\vrule height6pt width6pt depth0pt}

\newcommand{\Astar}{A^*}

\newcommand{\qDA}{{q^D_A}}
\newcommand{\qDB}{{q^D_B}}

\newcommand{\qA}{{q_A}}
\newcommand{\qB}{{q_B}}
\newcommand{\qC}{{q_C}}

\newcommand{\Aplus}{~\ {+}_{A}\ ~}
\newcommand{\Aminus}{~\ {-}_{A}\ ~}
\newcommand{\Atimes}{\ {\cdot}_{A}\ }

\newcommand{\zeroDA}{\tilde{0}^A}

\newcommand{\CStarplus}{~\ {+}_{C^*}\ ~}
\newcommand{\CStarminus}{~\ {-}_{C^*}\ ~}
\newcommand{\CStartimes}{\ {\cdot}_{C^*}\ }

\newcommand{\Bplus}{~\ {+}_{B}\ ~}
\newcommand{\Bminus}{~\ {-}_{B}\ ~}
\newcommand{\Btimes}{\ {\cdot}_{B}\ }
\newcommand{\zeroDB}{\tilde{0}^B}

\newcommand{\DDualA}{{D^{\star}A}}
\newcommand{\DDualB}{{D^{\star}B}}

\newcommand{\qDDualAA}{{q^{\star A}_A}}
\newcommand{\qDDualACStar}{{q^{\star A}_{C^*}}}

\newcommand{\qDDualBB}{{q^{\star B}_B}}
\newcommand{\qDDualBCStar}{{q^{\star B}_{C^*}}}

\newcommand{\CStar}{C^*}
\newcommand{\qCStar}{q_{C^*}}

\newcommand{\AStar}{A^*}
\newcommand{\qAStar}{q_{A^*}}

\newcommand{\BStar}{B^*}
\newcommand{\qBStar}{q_{B^*}}

\newcommand{\zeroDDualAA}{{\tilde{0}^{\star A}}}

\newcommand{\zeroDDualACStar}{\tilde{0}^{\star A,\CStar}}

\newcommand{\contract}{\lrcorner\,}

\newcommand{\conpairing}[1]{\left\{ #1\right\}}

\newcommand{\jet}{\mathfrak{J}}
\newcommand{\dev}{\mathfrak{D}}
\newcommand{\spray}[1]{\frkX(#1)}

\newcommand{\TE}{TE}

\newcommand{\EStar}{{E^*}}
\newcommand{\TStarE}{T^*E}

\newcommand{\TEStar}{TE^*}
\newcommand{\pE}{p_E}

\newcommand{\pEStar}{p_{E^*}}
\newcommand{\Tq}{T q }
\newcommand{\TqStar}{T q_* }
\newcommand{\rE}{r_E}
\newcommand{\cE}{c_E}
\newcommand{\rEStar}{r_{E^*}}
\newcommand{\cEStar}{c_{E^*}}
\newcommand{\TStarEStar}{T^*{E^*}}

\newcommand{\TStarM}{T^*M}
\newcommand{\pStar}{p^*}
\newcommand{\qStar}{q_*}

\newcommand{\TMplus}{~\ +_{\scriptstyle TM}\ ~}
\newcommand{\TMminus}{~\ -_{\scriptstyle TM}\ ~}
\newcommand{\TMtimes}{\ \cdot_{\scriptstyle TM}\ }

\newcommand{\qOmega}{{q_{\Omega}}}
\newcommand{\combineFunctor}{\frkC}

\newcommand{\doubleFunctor}{\frkD}
\newcommand{\idFunctor}{\mathbf{I}}

\newcommand{\Hom}{\mathrm{Hom}}
\newcommand{\gl}{\mathrm{gl}}

\newcommand{\splinear}{\mathrm{sl}}
\newcommand{\DVB}{\bf\mathrm{DVB}}

\newcommand{\UStar}{U^*}
\newcommand{\VStar}{V^*}
\newcommand{\KStar}{K^*}
\newcommand{\Id}{\mathrm{Id}}

\newcommand{\id}{I_{\dev}}
\newcommand{\jd}{J_{\dev}}
\newcommand{\e}{\mathbbm{e}}
\newcommand{\p}{\mathbbm{p}}

\begin{document}
\title{
{On Double Vector Bundles
\thanks
 {
 Research partially supported by NSFC (10871007,11026046,11001146) and Doctoral Fund. of MEC (20100061120096).
 }
} }
\author{ Zhuo Chen\\
Department of Mathematics, \\\vspace{2mm}Tsinghua University,
Beijing 100084, China\\ Zhangju Liu\\ Department of Mathematics and
LMAM\\\vspace{2mm}
Peking University, Beijing 100871, China\\ Yunhe Sheng\\
Department of Mathematics \\ Jilin University,
 Changchun 130012, Jilin, China\\
          { email: shengyh@jlu.edu.cn} }

\date{}


\maketitle

\begin{abstract}
 In this paper, we construct a category of short exact sequences of
 vector bundles and prove that it is equivalent to the
category of double vector bundles. Moreover, operations on double
vector bundles can be transferred to operations on the corresponding
short exact sequences. In particular, we study the duality theory of
double vector bundles in term of the corresponding short exact
sequences. Examples including the jet bundle and the Atiyah
algebroid are discussed.
\end{abstract}

\section{Introduction}

 The notion of a double vector bundle was
introduced by Pradines in \cite{Pradines:1968} and further studied
by Mackenzie in \cite{Mkz:1992,Mkz:2005}, Konieczna and
Urba$\rm\acute{n}$ski in \cite{KonieUrb:1989}, Grabowski in
\cite{grabowski} and Gracia-Saz and Mehta in \cite{alfonso} (this
list is not complete). As explained in \cite{Mkz:2005}, a double
vector bundle is essentially a vector bundle object in the category
of vector bundles. The most fundamental example  of double vector
bundles is the tangent   double vector bundle  $TE$  of a vector
bundle $E$, which is also known as the tangent double and had been
used in the connection theory since the late 1950s. It turns out
that the framework of double vector bundles is very convenient for
many important constructions like linear $1$-forms, linear Poisson
structures, special symplectic manifolds, linear connections and so
on (see
\cite{Besse:1978,Mkz:1992,Pradines:1974,YanoIshihara,KonieUrb:1989}
for more details).

The theory of duality of double vector bundles is decisively
different from that of ordinary vector bundles. A double vector
bundle has two duals which are themselves in duality. In fact, for a
double vector bundle $(D;A,B;M)_{C}$ (see Diagram
\eqref{Fig:DoubleVBDABC}), one can take the dual over $A$ and obtain
the double vector bundle $(\DDualA;A,\CStar;M)_{\BStar}$ shown as in
Diagram \eqref{Fig:DDualA}. In the mean time, one may also take the
dual over $B$ and obtain the double vector bundle
$(\DDualB;\CStar,B;M)_{\AStar}$ shown as in Diagram
\eqref{Fig:DDualB}. The observation that there is a natural duality
between $\DDualA$ and $\DDualB$ over $C^*$ is  illustrated by
\cite[Theorem 9.2.2]{Mkz:GTGA}.

The present paper is concerned with the algebraic feature of double
vector bundles. The main object introduced in this paper is  a
double vector bundle  sequence (DVB sequence), namely a short exact
sequence of vector bundles:
$$0\lon C\stackrel{\frke}{\lon}\Omega\stackrel{\frkp}{\lon} A\otimes
B\lon0 .$$ The main result of this paper is that such a DVB sequence
encapsulates the structure of a double vector bundle. In fact,  from
a DVB sequence, one can construct a double vector bundle
$(\doubleFunctor(\Omega);A,B;M)_C$ (this process is called the
double realization). Conversely, from a double vector bundle
$(D;A,B;M)_{C}$, we can also construct a DVB sequence $0\lon
C\lon\mathfrak{C}(D)\stackrel{\frkp}{\lon} A\otimes B\lon0$, which
is called the associated DVB sequence. Moreover, we show that the
category of double vector bundles (which will be denoted by ${\bf
DVB}$) and the category of DVB sequences (which will be denoted by
${\bf DVBS}$) are equivalent.

The dual form of a DVB sequence is called a DVB$^*$ sequence. In an
analogue way to the theory of duality of double vector bundles, we
introduce a method of taking duals to DVB$^*$ sequences. The dual
sequence obtained here (Sequence \eqref{seq:temp}) was also given by
Gracia-Saz and Mehta in \cite{alfonso}. It plays an important role
when the authors studied Lie algebroid structures on double vector
bundles. Furthermore, it turns out that a DVB$^*$ sequence has two
duals which are themselves in duality. Note that the pairing of such
dual DVB$^*$ sequences takes values in a vector bundle. Actually,
this idea already appeared in our work \cite{CLomni,CLS}, to study
certain geometric structures.

This paper is organized as follows. In Section 2 we give a review of
double vector bundles. In Section 3 we study DVB sequences of vector
bundles. In Subsection 3.1 we prove that one can construct a double
vector bundle associated to any DVB sequence (Theorem
\ref{thm:construction}). In Subsection 3.2 we prove that one can
also construct a DVB sequence associated to a double vector bundle.
Moreover, the category of double vector bundles and the category of
DVB sequences are equivalent (Theorem \ref{Thm:DVB1-1CM}). In
Subsection 3.3 we construct another short exact sequence of vector
bundles associated to a double vector bundle and prove that it is
isomorphic to the dual of the associated DVB sequence. In Section 4
we study the duality theory. In Subsection 4.1 we review the duality
theory of double vector bundles. In Subsection 4.2 we introduce the
notion of DVB$^*$ sequences and study their duals (taking values in
a vector bundle). We prove that after taking three steps of duals,
one can go back to the starting point (Theorem
\ref{Thm:TransLDual}). This is totally different from the usual
duality theory. In Subsection 4.3 we apply this method to the
associated DVB$^*$ sequence of a double vector bundle and prove that
this is the right process corresponding to taking the usual dual of
a double vector bundle (Theorem \ref{Thm:ADualCombine}). In Section
5 we give several examples including the Atiyah algebroid and the
jet bundle.



\section{ Double vector bundles}
Throughout the paper, we follow the conventions of \cite{Mkz:GTGA}.
\begin{defi}  A double vector bundle $(D;A,B;M)$ is a system of four
vector bundle structures
$$
\xymatrix{
  D \ar[d]_{\qDA} \ar[r]^{\qDB}
                & B \ar[d]^{\qB} &   \\
  A  \ar[r]^{\qA}
                & M         }
$$
in which $D$ has two vector bundles structures, on bases $A$ and
$B$, which are themselves vector bundles on $M$, such that each of
the four structure maps of each  vector bundle structure on $D$
(namely the bundle projection, zero section, addition and scalar
multiplication) is a morphism of vector bundles with respect to the
other structures.
\end{defi}

In the above figure, we refer to $A$ and $B$ as the side bundles of
$D$, and to $M$ as the double base. In the two side bundles, the
addition, scalar multiplication and subtraction are denoted by the
usual symbols $+$, juxtaposition, and $-$. We distinguish the two
zero-sections, writing $0^A: M\lon A$, $m\mapsto 0^A_m$, and $0^B:
M\lon B$, $m\mapsto 0^B_m$. We denote an element $d\in D$ by
$(d;a,b;m)$ to indicate that $\qDA(d)=a$, $\qDB(d)=b$,
$m=\qB(b)=\qA(a)$. We call
$D_m=(\qDA)\inverse(A_m)=(\qDB)\inverse(B_m)$ the \emph{slice} of
$D$ at $m$.

In the vertical bundle structure on $D$ with base $A$, the vector
bundle operations are denoted by $\Aplus$, $\Aminus$, $\Atimes$,
with $\zeroDA: A\lon D$, $a\mapsto \zeroDA_a$, for the zero-section.
Similarly, in the horizontal bundle structure on $D$ with base $B$
we write $\Bplus$, $\Bminus$, $\Btimes$, with $\zeroDB: B\lon D$,
$b\mapsto \zeroDB_b$.

The two structures on $D$, namely $(D,\qDB,B)$ and $(D,\qDA,A)$ will
also be denoted, respectively, by $\tilde{D}_B$ and $\tilde{D}_A$,  and called the
horizontal bundle structure and the vertical bundle structure.

The condition that each operation in $D$ is a morphism with respect
to the other is equivalent to the following equalities, known as the
\emph{interchange laws}.
\begin{eqnarray*}\label{Eqt:interchange1}
(d_1\Bplus d_2)\Aplus (d_3\Bplus d_4)&=& (d_1\Aplus d_3)\Bplus
(d_2\Aplus d_4),\\\label{Eqt:interchange2} t\Atimes(d_1\Bplus
d_2)&=&t\Atimes d_1\Bplus t\Atimes d_2,\\\label{Eqt:interchange3}
t\Btimes(d_1\Aplus d_2)&=&t\Btimes d_1\Aplus t\Btimes
d_2,\\\label{Eqt:interchange4} t\Atimes(s\Btimes d)&=&
s\Btimes(t\Atimes d),\\\label{Eqt:interchange5}
\zeroDA_{a_1+a_2}&=&\zeroDA_{a_1}\Bplus \zeroDA_{a_2},\\
\zeroDA_{ta}&=&t\Btimes \zeroDA_{a},\\\label{Eqt:interchange6}
\zeroDB_{b_1+b_2}&=&\zeroDB_{b_1}\Aplus \zeroDA_{b_2},\\
\zeroDB_{tb}&=&t\Atimes \zeroDB_{b}.
\end{eqnarray*}

We denote by $C$ the intersection of the two kernels:
$$
C=\set{c\in D~|~ \exists~ m\in M \mbox{ such that } \qDB(c)=0^B_m,
\quad \qDA(c)=0^A_m},
$$
which is called the core, and together with the map $\qC: c\mapsto
m$, $(C,\qC,M)$ is also a vector bundle over $M$. Although $C$ is
affiliated with $D$, in this paper we prefer to use
Diagram  (\ref{Fig:DoubleVBDABC}) below to emphasis the core of the
relevant double vector bundle.
\begin{equation}\label{Fig:DoubleVBDABC}
\xymatrix{
  D \ar[d]_{\qDA} \ar[r]^{\qDB}
                & B \ar[d]^{\qB} &   \\
  A  \ar[r]_{\qA}
                & M   & \ar[l]^{\qC} C  .       }
\end{equation}
A double vector bundle as above will also be  written as $(D;A,B;M)_{C}$.

An element $c\in C$ will be distinguished from its image in $D$,
which will be denoted by $\overline{c}$. The linear structures of
$C$ are inherited from those of $D$, due to the following facts:
\begin{equation}\label{Eqt:CoreProperty}
\overline{c+c'}=\overline{c}\Aplus\overline{c'}=\overline{c}\Bplus\overline{c'},
\quad \overline{rc}=r\Atimes\overline{c}=r\Btimes\overline{c}.
\end{equation}

Let $(D;A,B;M)_C$ be a double vector bundle given as above. The
\emph{flip} of $(D;A,B;M)_C$ is the double vector bundle
\begin{equation}\label{Fig:DoubleVBDBAC}
\xymatrix{
  D  \ar[d]_{\qDB} \ar[r]^{\qDA}
                & A \ar[d]^{\qA} &   \\
  B  \ar[r]_{\qB}
                & M   & \ar[l]^{\qC} C          }
\end{equation}
obtained by simply reversing the two side bundles. In this paper, we treat a double vector bundle and its flip as the same object.

\begin{defi}
\label{Def:morphismDouble} A morphism of double vector bundles
{\rm(see Diagram  (\ref{morphism}))}
    $$ (\varphi;f_A,f_B;f_M):\qquad(D;A,B;M)_C
    \rightarrow(D';A',B';M')_{C'}$$
consists of maps $\varphi$: $D\rightarrow D'$, $f_A:A\rightarrow
A'$, $f_B:B\rightarrow B'$, $f_M:M\rightarrow M'$, such that each of
$(\varphi,f_B)$, $(\varphi,f_A)$, $(f_A,f_M)$ and $(f_B,f_M)$ is a
morphism of the relevant vector bundles.
\begin{equation}\label{morphism}
\xymatrix@!0{
   D \ar[rr] \ar[dd] \ar[dr]_{\varphi}
      &  & B \ar'[d][dd]  \ar[dr]^{f_B}  &    \\
  & D' \ar[rr]\ar[dd]
      &  & B' \ar[dd] & \\
   A \ar'[r][rr]\ar[dr]_{f_A}
      &  & M\ar[dr]^{f_M}               \\
  & A' \ar[rr]
      &  & M'    .     }
\end{equation}
\end{defi}
In this definition,   $(f_C,f_M)$, where $f_C=\varphi|_{C}:C\rightarrow
C'$, is also a morphism of the associated cores.

We denote the category of double vector bundles by ${\bf DVB}$.

\begin{ex}\label{Ex:trivial1}\rm Consider $D=A\times_M B\times_M C$,
the trivial double vector bundle over $A$ and $B$ with core $C$ (see
\cite[Example 9.1.4]{Mkz:GTGA}). The bundle projections $\qDA$ and
$\qDB$ are exactly the projections to $A$ and $B$ respectively. The
linear structures are given by
$$
r\Atimes(a,b_1,c_1)\Aplus (a,b_2,c_2)=(a,rb_1+b_2,rc_1+c_2);
$$
$$
r\Btimes(a_1,b,c_1)\Bplus (a_2,b,c_2)=(ra_1+a_2,b,rc_1+c_2).
$$

\end{ex}

In \cite{grabowski} and \cite{alfonso}, it is proved that
any double vector bundle is isomorphic to a trivial double vector
bundle.

\section{DVB sequences}

In this section, $\qOmega: \Omega\lon M$, $\qA: A\lon M$, $\qB:
B\lon M$ and $\qC:C\lon M$ are vector bundles over $M$.
\begin{defi}
A DVB sequence\footnote{DVB is the abbreviation of ``double vector
bundle''.}   on $M$ is an exact sequence of vector bundles:
\begin{equation}\label{CoreCombination}
\xymatrix@C=0.5cm{
  0 \ar[r] & C \ar[rr]^{\frke} && \Omega \ar[rr]^{\frkp} && A\otimes B \ar[r] & 0 }
\end{equation} over the identity map $\Id_{M}$.
\end{defi}

Such a DVB sequence will be denoted by
$(\Omega\stackrel{\frkp}{\longrightarrow} A\otimes B;M)_C$. We refer to the vector bundles $A$, $B$ as side bundles, and $C$ the core.

\begin{defi}\label{Def:morphismCombination}
A morphism of DVB sequences
$$(\varpi;f_A,f_B;f_M):\quad
(\Omega\stackrel{\frkp}{\longrightarrow} A\otimes B;M)_C\lon
(\Omega'\stackrel{\frkp'}{\longrightarrow} A'\otimes B';M')_{C'},
$$
consists of maps $\varpi: \Omega'\lon \Omega$,  $f_A: A\lon A'$,
$f_B: B\lon B'$ and  $f_M: M\lon M'$, each of $(\varpi,f_M)$,
$(f_A,f_M)$, $(f_B,f_M)$
 is a morphism of the relevant vector bundles, such that the diagram
\begin{equation}\label{MorphismCombination}
\xymatrix@!0{
   \Omega \ar[rr]|-{\frkp} \ar[dd] \ar[dr]_{\varpi}
      &  & A\otimes B \ar'[d][dd]  \ar[dr]^{f_A\otimes f_B}  &    \\
  & \Omega' \ar[rr]|-{\frkp'}\ar[dd]
      &  & A'\otimes B' \ar[dd] & \\
   M \ar@{=}'[r][rr]\ar[dr]_{f_M}
      &  & M\ar[dr]^{f_M}               \\
  & M' \ar@{=}[rr]
      &  & M'        }
\end{equation}
commutes.
\end{defi}
In this definition, for $f_C=\varpi|_{C}:C\rightarrow C'$,
$(f_C,f_M)$ is also a morphism of vector bundles.

It is obvious that DVB sequences together with the above morphisms form a category, which will be denoted by ${\bf DVBS}$.

\subsection{From DVB sequences to double vector bundles}

We show how to construct a double vector bundle out of a DVB sequence $(\Omega\stackrel{\frkp}{\longrightarrow} A\otimes B;M)_C$.
\begin{thm}\label{thm:construction}
Given a DVB sequence
$(\Omega\stackrel{\frkp}{\longrightarrow} A\otimes B;M)_C$, there is an associated
 double vector bundle
\begin{equation}\label{DoubleRealization}
 \xymatrix{
  \doubleFunctor(\Omega) \ar[d]_{\qDA} \ar[r]^{\qDB}
                & B \ar[d]^{\qB  } &   \\
  A  \ar[r]^{\qA}
                & M.        }
\end{equation}
Here $\doubleFunctor(\Omega)$ is given by
$$
\doubleFunctor(\Omega)\defbe \set{(\omega,a,b)\in \Omega\times_M
A\times_M B| \frkp(\omega)=a\otimes b}.
$$
The vertical bundle projection $\qDA: D\lon A$ is given by
$$
\qDA: (\omega,a,b)\mapsto a,\quad \forall (\omega,a,b)\in
\doubleFunctor(\Omega).
$$
The horizontal bundle projection $\qDB: D\lon B$ is given by
$$
\qDB: (\omega,a,b)\mapsto b,\quad \forall (\omega,a,b)\in
\doubleFunctor(\Omega).
$$
\end{thm}
\pf Since $\frkp$ is surjective, both $\qDA$ and $\qDB$ are
surjections. For the vector bundle structure  over $A$, we define (for any
fixed $a\in A$)
$$
r\Atimes (\omega_1,a,b_1)\Aplus (\omega_2,a,b_2)\defbe
(r\omega_1+\omega_2,a,rb_1+b_2).
$$
Similarly, the vector bundle structure over $B$ are defined by  (for any
fixed $b\in B$):
$$
r\Btimes (\omega_1,a_1,b)\Bplus (\omega_2,a_2,b)\defbe
(r\omega_1+\omega_2,ra_1+a_2,b).
$$
For $\doubleFunctor(\Omega)$, the zero element above
$a\in A$ is $\zeroDA_a=(0,a,0)$ and the zero element above $b\in B$ is
$\zeroDB_b=(0,0,b)$. The core of $(\doubleFunctor(\Omega);A,B;M)$ can be identified with the core $C=Ker(\frkp)$ of the DVB sequence
$(\Omega\stackrel{\frkp}{\longrightarrow} A\otimes B;M)_C$. In fact,
any $c\in C$ is embedded into $\doubleFunctor(\Omega)$ by
$\overline{c}=(\frke(c),0,0)$. It is routine to check that the
interchange laws are satisfied.  \qed

\begin{defi} We call the double vector bundle $(\doubleFunctor(\Omega);A,B;M)_C$ shown as in Diagram  {\rm(\ref{DoubleRealization})} the double
realization of the DVB sequence
$(\Omega\stackrel{\frkp}{\longrightarrow} A\otimes B;M)_C$.
\end{defi}

For any morphism of DVB sequences
$$(\varpi;f_A,f_B;f_M):~~(\Omega\stackrel{\frkp}{\longrightarrow}
A\otimes B;M)_C\longrightarrow
(\Omega'\stackrel{\frkp}{\longrightarrow} A'\otimes B';M)_{C'},$$
define $\doubleFunctor(\varpi):\doubleFunctor(\Omega)\longrightarrow
\doubleFunctor(\Omega')$ by
$$
\doubleFunctor(\varpi) (\omega,a,b)=(\varpi(\omega),
f_A(a),f_B(b)),\quad\forall~ (\omega,a,b)\in \doubleFunctor(\Omega).
$$
It is easily seen that
$$
(\doubleFunctor(\varpi),f_A,f_B;f_M):(\doubleFunctor(\Omega);A,B;M)_C\longrightarrow
(\doubleFunctor(\Omega');A',B';M')_{C'}
$$
is a morphism of double vector bundles.\vspace{2mm}

\begin{lem}\label{lem:Dfunctor}
With the above notations, \begin{eqnarray*}\doubleFunctor:&
(\Omega\stackrel{\frkp}{\longrightarrow} A\otimes B;M)_C
\dashrightarrow (\doubleFunctor(\Omega);A,B;M)_C,\\
\doubleFunctor:&(\varpi;f_A,f_B;f_M)\dashrightarrow
(\doubleFunctor(\varpi),f_A,f_B;f_M)
\end{eqnarray*}
is a covariant functor {\rm (called the doubling functor)} from the
category of DVB sequences ${\bf DVBS}$ to the category of double
vector bundles ${\bf DVB}$.
\end{lem}
The proof is straightforward and thus omitted.

\subsection{From double vector bundles to DVB sequences}

In this part we show that a double vector bundle also
yields a DVB sequence. Suppose that we are given a double vector
bundle $(D;A,B;M)_{C}$ as in Diagram  (\ref{Fig:DoubleVBDABC}). Fix
$m\in M$ and let $\widetilde{D}_m$ denote the free vector space
generated by all the elements in the slice $D_m$. In other words, an
element of $\widetilde{D}_m$ is a formal linear combination of some
elements of $D_m$ with coefficients in $\Real$, i.e. it has an
expression of the form
$$\Sigma_{i=1}^l r_i\cdot d_i, \quad r_i\in
\Real, d_i\in D_m\,.
$$
Let $\overline{D}_m$ be the subspace of $\widetilde{D}_m$ generated
by elements of the following types:
\begin{equation}\label{Eqv1}
 \left\{\begin{array}{l}d_1 + d_2 -
(d_1\Aplus d_2), \\ d_3 + d_4 -(d_3\Bplus d_4),\\
r\cdot d - (r\Atimes d), \\ r\cdot d - (r\Btimes
d).\end{array}\right.
\end{equation}
Let $\combineFunctor(D)_m= \widetilde{D}_m/\overline{D}_m$. For each
$d\in D_m$, let $[d]$ denote its image in $\combineFunctor(D)_m$.
Then, $\combineFunctor(D)_m$ is generated by elements of the form
$[d]$, and by definition, one has
$$
[d_1\Aplus d_2]=[d_1]+[d_2],\quad [d_3\Bplus d_4]=[d_3]+[d_4],\quad
[r\Atimes d]=[r\Btimes d]=r[d].
$$
In particular, we know that
$$
[\zeroDA_a]=[\zeroDB_b]=0,\quad\forall~ a\in A_m,\ b\in B_m\,.
$$

\begin{pro}\label{Pro:dimOmegaD}The vector space
$\combineFunctor(D)_m$ has dimension $\dim A \times
\dim B + \dim C$.
\end{pro}
\pf Since any double vector bundle can be trivialized, it suffices
to assume $D_m= A_m\times B_m\times C_m$. In this case, the subspace
$\overline{D}_m$ of $\widetilde{D}_m$ is generated by elements of
the form (according to \eqref{Eqv1})
$$
(a,b_1,c_1)+(a,b_2,c_2)-(a,b_1+b_2,c_1+c_2),
$$
$$
(a_3,b, c_3)+(a_4,b,c_4)-(a_3+a_4,b,c_3+c_4),
$$
$$
r\cdot (a,b,c)-(a,rb,rc),
$$
$$
r\cdot (a,b,c)-(ra,b,rc).
$$
Hence the quotient space $ \combineFunctor(D)_m=
\widetilde{D}_m/\overline{D}_m$ is isomorphic to $A_m\otimes B_m
\oplus C_m$ (c.f. \cite[Chp.2]{Atiyah}) by canonically identifying $
[(a,b,c)]$ with $ a\otimes b + c$.
This completes the proof. \qed\vspace{3mm}




It is also clear that the map $[d]\mapsto a\otimes b$, for any $
(d;a,b,m)\in D_m$, extends to a linear map
\begin{equation}\label{frkpm}
\frkp_m: \combineFunctor(D)_m \lon A_m\otimes B_m\,,
\end{equation}
which is surjective.

In summary, we have
\begin{thm} The  set $\combineFunctor(D)=\cup_{m\in M}\combineFunctor(D)_m$ is a vector bundle
over $M$ and fit into the following DVB sequence:
\begin{equation}\label{CoreCombinationDABC}
\xymatrix@C=0.5cm{
  0 \ar[r] & C \ar[rr]^{\frke} && \combineFunctor(D) \ar[rr]^{\frkp} && A\otimes B \ar[r] & 0
  }.
\end{equation}
where $\frkp$ is induced by $ [d]\mapsto \qDA(d)\otimes
\qDB(d)$ and $\frke$ is given by $c\mapsto [\overline{c}]$.
\end{thm}

 We call
{\rm(\ref{CoreCombinationDABC})} the  {\bf
associated DVB sequence} of the double vector bundle $(D;A,B;M)_{C}$.

Let
$$ (\varphi;f_A,f_B;f_M):\quad(D;A,B;M)_C
    \rightarrow(D';A',B';M')_{C'}
$$
be a morphism of double vector bundles, define
$\combineFunctor(\varphi): \combineFunctor(D)\lon
\combineFunctor(D')$ by setting
$$
\combineFunctor(\varphi)([d])= [\varphi(d)].
$$
Obviously, it is a morphism of vector bundles and satisfies
$$
\frkp' \circ \combineFunctor(\varphi)= (f_A\otimes f_B)\circ \frkp.
$$
Hence $(\combineFunctor(\varphi);f_A,f_B;f_M)$ is a morphism of DVB
sequences.\vspace{2mm}

\begin{lem}\label{lem:Cfunctor}
With the above notations,
\begin{eqnarray*}
\combineFunctor: & (D;A,B;M)_C \dashrightarrow
 (\combineFunctor(D)\stackrel{\frkp}{\longrightarrow} A\otimes
 B;M)_C,\\
\combineFunctor:&(\varphi;f_A,f_B;f_M)\dashrightarrow
(\combineFunctor(\varphi);f_A,f_B;f_M)
\end{eqnarray*}
is a covariant functor  from the category of double vector bundles
${\bf DVB}$ to the category of DVB sequences ${\bf DVBS}$.
\end{lem}
Again we omit the proof of this lemma. The main result of this  paper is the following theorem.
\begin{thm}\label{Thm:DVB1-1CM}
The category of double vector bundles ${\bf DVB}$ and the category
of DVB sequences ${\bf DVBS}$ are equivalent.
\end{thm}
\pf For the preceding $\doubleFunctor$-functor and the
$\combineFunctor$-functor, respectively, given by Lemma \ref{lem:Dfunctor} and
\ref{lem:Cfunctor}, we   establish a natural transformation $t$
from the identity functor $\idFunctor_{\DVB}$ to
$\doubleFunctor\combineFunctor$. For any double vector bundle
$(D;A,B;M)$ and $(d;a,b;M)\in D$, as a canonical element $[d]\in
\combineFunctor(D)$, $([d],a,b)\in
\doubleFunctor\combineFunctor(D)$, define
$t_D:D\longrightarrow\doubleFunctor\combineFunctor(D)$ by setting
$$
t_D( d)= ([d],a,b).
$$
Since it preserves the side bundles $A$, $B$ and the core $C$,
$$
(t_D;\Id_A,\Id_B;\Id_M): ~~(D;A,B;M)_C\ \lon\
(\doubleFunctor\combineFunctor(D);A,B;M)_C
$$
 must be an isomorphism. Further, $t$ is natural in the  sense
that the following diagram
$$
\begin{CD}
(D;A,B;M)_C @>{t_D}>> (\doubleFunctor\combineFunctor(D);A,B;M)_C \\
@V{(\varphi;f_A,f_B;f_M)}VV @VV({(\doubleFunctor\combineFunctor)}(\varphi);f_A,f_B;f_M)V \\
(D';A',B';M')_{C'} @>{t_{D'}}>>
(\doubleFunctor\combineFunctor(D');A',B';M')_{C'}
\end{CD}
$$
is commutative, for every morphism of double vector bundles
${(\varphi;f_A,f_B;f_M)}$.

On the other hand, there is  a natural transformation $\pi$
from $\combineFunctor\doubleFunctor$ to the identity functor
$\idFunctor_{\bf DVBS}$. In fact, for any DVB sequence
$(\Omega\stackrel{\frkp}{\longrightarrow} A\otimes B;M)_C$, we can
establish a standard isomorphism
$$(\pi_{\Omega};\Id_A,\Id_B;\Id_M):\quad
(\combineFunctor
\doubleFunctor(\Omega)\stackrel{\frkp}{\longrightarrow} A\otimes
B;M)_C\ \lon\ (\Omega\stackrel{\frkp}{\longrightarrow} A\otimes
B;M)_C.
$$
This is done as follows. Write $D=\doubleFunctor(\Omega)$ and recall
that an element in $D_m$ is of the form
$$
d=(\omega,a,b),\quad \mbox{ satisfying }\ \frkp(\omega)=a\otimes b,\
\quad\mbox{ where }\ \omega\in \Omega_m,a\in A_m,b\in B_m.
$$
So we define a map $\widetilde{\pi}$: $\widetilde{D}_m\lon \Omega$
(recall that $\widetilde{D}_m$ is freely generated by $D_m$) by
setting $(\omega,a,b)\mapsto \omega$ and then it linearly extends to
$\widetilde{D}_m$. It is easy to check that $\widetilde{\pi}$ sends
all the four types of elements (\ref{Eqv1})  to zero. Therefore,
$\widetilde{\pi}$ induces a well defined linear map $\pi_{\Omega}$
from $\combineFunctor(D)_m= \widetilde{D}_m/\overline{D}_m$ to
$\Omega$, i.e.
$$
\pi_\Omega([(\omega,a,b)])=\omega.
$$
Obviously, $\pi_{\Omega}:\combineFunctor \doubleFunctor(\Omega)\lon
\Omega$ is a morphism of vector bundles.

One can directly verify that the diagram
$$
\xymatrix{
 0 \ar[r] & C \ar[d]_{\Id_C}  \ar[rr]^{\frke} &&
                \combineFunctor \doubleFunctor(\Omega) \ar[d]_{\pi_{\Omega}}  \ar[rr]^{\frkp} && A\otimes B \ar[d]_{\Id}  \ar[r]  & 0   \\
 0 \ar[r] & C  \ar[rr]^{\frke} &&
                \Omega  \ar[rr]^{\frkp} && A\otimes B \ar[r]  & 0,
                }
$$
commutes and it implies that $\pi_\Omega$ is an isomorphism of
vector bundles. Finally, $\pi$ is natural because for any morphism
$(\varpi;f_A,f_B;f_M)$ of DVB sequences, we have the following
commutative diagram:
$$
\begin{CD}
(\combineFunctor \doubleFunctor(\Omega)\stackrel{\frkp}{\longrightarrow} A\otimes B;M)_C @>{\pi_{\Omega}}>> (\Omega\stackrel{\frkp}{\longrightarrow} A\otimes B;M)_C \\
@V{((\combineFunctor
\doubleFunctor)(\varpi);f_A,f_B;f_M)}VV @VV{(\varpi;f_A,f_B;f_M)}V \\
(\combineFunctor
\doubleFunctor(\Omega')\stackrel{\frkp'}{\longrightarrow} A'\otimes
B';M')_{C'} @>{\pi_{\Omega'}}>>
(\Omega'\stackrel{\frkp'}{\longrightarrow} A'\otimes B';M')_{C'}.
\end{CD}
$$
This finishes the proof of the equivalence of the categories ${\bf
DVB}$ and   ${\bf DVBS}$. \qed

\subsection{The dual of the associated DVB sequence}
Consider the dual of the associated  DVB sequence
(\ref{CoreCombinationDABC}):
\begin{equation}\label{DualCoreCombination}
\xymatrix@C=0.5cm{
  0 \ar[r] & \AStar\otimes \BStar \ar[rr]^{\frkp^*} &&
\combineFunctor(D)^*
   \ar[rr]^{\frke^*} && \CStar \ar[r] & 0 }.
\end{equation}
In this part of the paper, we give an explanation of this sequence.

\begin{defi}
For the double vector bundle $(D;A,B;M)_C$ shown as in Diagram
{\rm(\ref{Fig:DoubleVBDABC})}, a function $\sigma$ on the slice
$D_m$ is said to be \textbf{double-linear}, if it satisfies the
following equalities
\begin{eqnarray*}\label{Eqt:altsigma1}
\sigma(d_1\Aplus d_2)&=&\sigma(d_1)+\sigma(d_2),\quad
    \mbox{ if } \qDA(d_1)=\qDA(d_2);\\\label{Eqt:altsigma2}
\sigma(r \Atimes d)&=&r \sigma(d);\\\label{Eqt:altsigma4}
\sigma(d_3\Bplus d_4)&=&\sigma(d_3)+\sigma(d_4),\quad\mbox{ if }
    \qDB(d_3)=\qDB(d_4);\\\label{Eqt:altsigma5}
\sigma(r \Btimes d)&=&r \sigma(d);
\end{eqnarray*} where $d,d_1,d_2,d_3,d_4\in D_m$, $r\in\Real$.
\end{defi}

For any $\theta\in\AStar_m\otimes\BStar_m$,  define a function
$i(\theta)$ on $D_m$ by
\begin{equation}\label{Eqt:AStarBStarFunctionOnSlice}
i(\theta)(d)\defbe \pairing{\theta\ ,\
\qDA\otimes\qDB(d)}=\pairing{\theta,a\otimes b},\quad \forall~
(d;a,b;m)\in D_m\,,
\end{equation}
which is obviously double-linear.

Let $(\spray{D})_m$ be the collection of double linear functions on
$D_m$. With respect to the usual addition and multiplication of
functions, $(\spray{D})_m$ is a vector space. It is easy to see that
$$\spray{D}=\bigcup_{m\in M} (\spray{D})_m$$ is a vector
bundle over $M$ and $i: \AStar\otimes\BStar\lon \spray{D}$ given by
(\ref{Eqt:AStarBStarFunctionOnSlice}) is a vector bundle morphism.

\begin{pro}\label{prop:XD}Any double-linear function $\sigma:D_m\lon \Real$
uniquely determines some $\chi\in \CStar_m$, such that
$$
\sigma(\overline{c})=\pairing{\chi,c}, \quad\forall~ c\in C_m\,.
$$
Write $\chi=j(\sigma)$, then $j$ is a vector bundle morphism and
\begin{equation}\label{SprayDABC}
\xymatrix@C=0.5cm{
  0 \ar[r] & \AStar\otimes \BStar \ar[rr]^{i} &&
\spray{D}
   \ar[rr]^{j} && \CStar \ar[r] & 0 }
\end{equation}
is an exact sequence.
\end{pro}

\pf By relation (\ref{Eqt:CoreProperty}), $\sigma|_{C_m}$ is a
linear function on $C_m$ and this clearly defines a linear map $j$:
$\sigma\lon \sigma|_{C_m}\in \CStar_m$. To see that $j$ is
surjective it is sufficient to work locally. In fact, if we assume
that $D_m= A_m\times B_m\times C_m$, then a double linear function
on $D_m$ is indeed an element in $\AStar\otimes \BStar\oplus
\CStar$. In this case, $j$ is the projection to $\CStar$. It is also
trivial that $i$ is an injection. \qed\vspace{3mm}

Below we establish a canonical way to identify the sequence
{\rm(\ref{SprayDABC})}   with the dual  of the associated DVB
sequence {\rm(\ref{CoreCombinationDABC})}. In fact, there is a
standard pairing between $\combineFunctor(D)$ and $\spray{D}$:
$$
\pairing{\cdot,\cdot}:\ \quad\combineFunctor(D)~~\times_M
\spray{D}~~\lon \Real,
$$
induced  by $$
\pairing{[d],\sigma}\defbe \sigma(d),
$$
for all $(d;a,b;m)\in D_m$ and $\sigma\in (\spray{D})_m$.
It is not hard to see that this well defines a pairing between
$\combineFunctor(D)$ and $\spray{D}$. Moreover, it
satisfies
\begin{eqnarray}\label{PairingSpecial1}
\left\{\begin{array}{ll}\pairing{\omega,i(\theta)}=\pairing{\frkp(\omega),\theta},\quad&
\forall~ \omega\in \combineFunctor(D)_m,\ \theta\in \AStar_m\otimes\BStar_m;\\
\pairing{\frke(c),\sigma}=\pairing{c,j(\sigma)},\quad&\forall~ c\in
C_m,\ \sigma\in(\spray{D})_m.\end{array}\right.
\end{eqnarray}

Fix $\sigma\in\frkX(D)_m$, the condition
``$\pairing{[d],\sigma}=\sigma(d)=0$, for all $d\in D_m$'' clearly
implies $\sigma=0$.   Thus the pairing  must be nondegenerate.

\section{Duality theory}

\subsection{Duality of double vector bundles}

Dualizing the vertical structure on $D$ leads again to a double
vector bundle $(\DDualA;A,\CStar;M)_{\BStar}$, called the vertical
dual or dual over $A$ of (\ref{Fig:DoubleVBDABC}), which is denoted
by
\begin{equation}\label{Fig:DDualA}
\xymatrix{
  \DDualA \ar[d]_{\qDDualAA} \ar[r]^{\qDDualACStar}
                & \CStar \ar[d]^{\qCStar} &   \\
  A  \ar[r]_{\qA}
                & M   & \ar[l]^{\qBStar} \BStar .        }
\end{equation}
The vertical structure in (\ref{Fig:DDualA}) is the usual dual of
the bundle structure on $\tilde{D}_A$, and $\qCStar$ is the usual
dual of $\qC: C\lon M$. The additions and scalar multiplications in
the side bundles will be denoted by $\Aplus$, $\Atimes$, $\Aminus$
and $\CStarplus$, $\CStartimes$, $\CStarminus$. The zero of
$\DDualA$ above $a\in A$ is denoted by $\zeroDDualAA_a$.

The vector bundle projection $\qDDualACStar$: $\DDualA\lon\CStar$ is
defined by
\begin{equation}\label{Eqt:qDDualACStar}
\pairing{\qDDualACStar(\Phi),c}=\pairing{\Phi,\zeroDA_a\Bplus\overline{c}},
\quad\forall~ a\in A_m, c\in C_m, \Phi\in (\qDDualAA)\inverse(a).
\end{equation}
The addition $\CStarplus$ in $\DDualA\lon\CStar$ is defined by
\begin{equation}\label{eqt:Cstarplus}
\pairing{\Phi \CStarplus \Phi', d\Bplus d' } =\pairing{\Phi,d}+
\pairing{\Phi',d'}.
\end{equation}
Here $\Phi\in(\qDDualAA)\inverse(a)$,
$\Phi'\in(\qDDualAA)\inverse(a')$, and $d\in(\qDA)\inverse(a)$,
$d'\in(\qDA)\inverse(a')$. Note the important hypothesis that
$\qDDualACStar(\Phi)=\qDDualACStar(\Phi')$, which ensures that
$\CStarplus$ is well-defined.

Similarly, we define
 $$
 \pairing{r\CStartimes \Phi,
d}=r\pairing{\Phi,\frac{1}{r}\Btimes d},
$$
for $r\neq 0$ and $d\in D$ with $\qDA(d)=r \cdot \qDDualAA(\Phi)$.

The zero above $\chi \in \CStar_m$ is denoted by
$\zeroDDualACStar_\chi$ and is defined by
$$
\pairing{\zeroDDualACStar_\chi,\zeroDB_b\Aplus \overline{c}}=
\pairing{\chi,c}, \quad \forall~ b\in B_m,c\in C_m.
$$
The core element $\overline{\psi}$ corresponding to $\psi\in
\BStar_m$ is
$$
\pairing{\overline{\psi},\zeroDB_b\Aplus \overline{c}}=
\pairing{\psi,b}, \quad\mbox{i.e.}\quad
\pairing{\zeroDDualAA_a\CStarplus
\overline{\psi},d}\defbe\pairing{\psi,\qDB(d)},\quad\forall~ d\in D.
$$

There is of course also a horizontal dual
$(\DDualB;\CStar,B;M)_{\AStar}$:
\begin{equation}\label{Fig:DDualB}
\xymatrix{
  \DDualB \ar[d]_{\qDDualBCStar} \ar[r]^{\qDDualBB}
                & B \ar[d]^{\qB} &   \\
  \CStar  \ar[r]_{\qCStar}
                & M   & \ar[l]^{\qAStar} \AStar ,        }
\end{equation}
defined in an analogous way.

\subsection{Duality of DVB$^*$ sequences}
In this part of the paper, $U$, $V$, $K$ and $\Pi$ are all vector
bundles over the base space $M$. The dual form of DVB sequences will
be called DVB$^*$ sequences.
\begin{defi}
An exact sequence of vector bundles over $M$ as follows
\begin{equation}\label{seq:prejet}
 \xymatrix@C=0.5cm{
  0 \ar[r] & U\otimes V \ar[rr]^{i} && \Pi \ar[rr]^{j} &&
  K \ar[r] & 0, }
\end{equation}
is called a DVB$^*$ sequence. Again, we refer to $U$ and $V$ as the
side bundles.
\end{defi}

Obviously, the short exact sequence  (\ref{SprayDABC}) is a DVB$^*$
sequence, which we call the \textbf{associated DVB$^*$ sequence} of
the double vector bundle $(D;A,B;M)_C$.
\begin{defi}
Let
\begin{equation}\label{seq:prejet1}
\xymatrix@C=0.5cm{
  0 \ar[r] & U\otimes \KStar \ar[rr]^{e} && \Delta\ar[rr]^{p} &&
  \VStar \ar[r] & 0  }
\end{equation}
be a DVB$^*$ sequence. The DVB$^*$ sequences (\ref{seq:prejet}) and
\eqref{seq:prejet1} are said to be in duality with respect to $U$,
or $U$-duality, if
 there is a bilinear pairing taking values in $U$:
$$
\conpairing{\cdot,\cdot}_{U}: \Delta\times_M \Pi\lon U,
$$
satisfying the following two equalities
\begin{eqnarray}\label{Eqt:PairingConjugate1}
\left\{\begin{array}{ll}\conpairing{\epsilon,i(\theta)}_{U}=
p(\epsilon) \contract \theta,\quad&\forall~
\epsilon\in\Delta,\theta\in U\otimes V\,,\\
\conpairing{e(\kappa),\sigma}_{U}=j(\sigma)\contract \kappa,
\quad&\forall~ \sigma\in\Pi,\kappa\in U\otimes
\KStar\,.\end{array}\right.
\end{eqnarray}
\end{defi}

\begin{rmk}\rm
 It is easy to check
the following facts:
$$
\conpairing{\cdot,\sigma}_{U}=0\quad\Longleftrightarrow\quad
\sigma=0;
$$
$$
\conpairing{\epsilon,\cdot}_{U}=0\quad\Longleftrightarrow\quad
\epsilon=0.
$$
Thus the $U$-valued pairing $\conpairing{\cdot,\cdot}_{U}$ in this
definition must be nondegenerate.
\end{rmk}

 Below we use a constructive approach to show the existence of such duals of a given DVB$^*$ sequence.    Again consider the  DVB$^*$ sequence
(\ref{seq:prejet}) and its dual form:
$$
\xymatrix@C=0.5cm{
  0 \ar[r] & \KStar \ar[rr]^{j^*} && \Pi^* \ar[rr]^{i^*} &&
  \UStar\otimes \VStar \ar[r] & 0. }
$$
Applying the ``$ U\otimes$'' operation to this sequence, we get
another exact sequence
$$
\xymatrix@C=0.5cm{
  0 \ar[r] & U\otimes \KStar \ar[rr]^{\Id_U\otimes j^*} && U\otimes\Pi^*
   \ar[rr]^{\Id_U\otimes i^*} &&
 U\otimes \UStar\otimes \VStar \ar[r] & 0. }
$$
Using the standard decomposition $U\otimes \UStar=\gl(U)\cong
\splinear(U)\oplus \Real \Id_U$, the right side of the above
sequence becomes $\splinear(U)\otimes \VStar \oplus \VStar$. So we
are able to obtain a sub-vector bundle of $U\otimes\Pi^*$, which is
the pull back of $\VStar$ and will be denoted by $\Pi^{*}_U$. Moreover, we
have the following exact sequence:
\begin{eqnarray}\label{Seq:prejet2}
&\xymatrix@C=0.5cm{
  0 \ar[r] & U\otimes \KStar \ar[rr]^{e}
  && \Pi^{*}_U \ar[rr]^{p} &&
  \VStar \ar[r] & 0, }&
\end{eqnarray}
where $e=\Id_U\otimes j^*,\; p=(\Id_U\otimes i^*)|_{\Pi^{*}_U}.$

Since  $\Pi^{*}_U$ is the pull back of $V^*$, it can be directly
described as follows:
$$
(\Pi^{*}_U)_m=\set{\epsilon\in \Hom(\Pi_m,U_m)\,|\,\exists~ v^*\in
\VStar_m, \mbox{ s.t., } \epsilon\circ
i(\theta)=v^*\contract\theta,~\forall~\theta\in U\otimes V}.
$$

It is a straightforward verification to prove the following fact:
\begin{pro}
The DVB$^*$ sequences \eqref{seq:prejet} and \eqref{Seq:prejet2} are
in $U$-duality.
\end{pro}

In a similar manner, one is able to get another DVB$^*$ sequence
\begin{eqnarray*}
&\xymatrix@C=0.5cm{
  0 \ar[r] &  \KStar \otimes V\ar[rr]^{e}
  && \Pi^{*}_V \ar[rr]^{p} &&
  U^* \ar[r] & 0,}
\end{eqnarray*}
which is a $V$-dual of \eqref{seq:prejet}.\vspace{2mm}


\begin{rmk}\rm In the special case that $U=M\times\Real$, the
DVB$^*$ sequence \eqref{seq:prejet} is of the form
$$
 \xymatrix@C=0.5cm{
  0 \ar[r] &  V \ar[rr]^{i} && \Pi \ar[rr]^{j} &&
  K \ar[r] & 0. }
$$
It is easy to see that its $(M\times\Real)$-dual
$\Pi^{*}_{M\times\Real}$ is exactly $\Pi^*$. So we get the usual
dual of an exact sequence.
\end{rmk}

Let us introduce another operation of DVB$^*$ sequences. In short,
transpositions  exchange the side bundles and add a minus sign. In
specific,  the {\bf transposition} of the DVB$^*$ sequence
(\ref{seq:prejet}) is the following DVB$^*$ sequence:
\begin{eqnarray}\label{seq:prejettrans}
\xymatrix@C=0.5cm{
  0 \ar[r] & V\otimes U \ar[rr]^{i^t} && \Pi \ar[rr]^{j} &&
  K \ar[r] & 0, }
\end{eqnarray}
where $i^t$ is given  by
$$
v\otimes u\ \ \mapsto\ \ -i(u\otimes v),\quad \forall ~u\in U,v\in
V.
$$

A DVB$^*$ sequence remains unchanged if we take twice
transpositions. An interesting phenomenon is the next theorem, which
claims that after taking \textbf{three} steps of  duals   with
respect to different side bundles, one gets the transposition.

\begin{thm}\label{Thm:TransLDual}
Let
\begin{equation}\label{seq:prejet3}
0\lon V\otimes \KStar
\stackrel{r}{\longrightarrow}(\Pi^*_U)^*_{\KStar}
\stackrel{q}{\longrightarrow} \UStar \lon 0
\end{equation}be the $K^*$-dual of \eqref{Seq:prejet2}.
Then the DVB$^*$ sequences \eqref{seq:prejet3} and
\eqref{seq:prejettrans} are in $V$-duality, or, in terms of a
 diagram:
\begin{eqnarray}\nonumber \fbox{\eqref{seq:prejet}} &\quad\stackrel{\mbox{\emph{ transposition}}}{
\Longleftrightarrow}\quad& \fbox{\eqref{seq:prejettrans}}\\\label{Dig:square}
\\\nonumber U \mbox{\emph{-dual}}{\bigg \Updownarrow} &
& {\bigg \Updownarrow}V
\mbox{\emph{-dual}}\\\nonumber\quad\\\nonumber
\fbox{\eqref{Seq:prejet2}} &\stackrel{\KStar\mbox{\emph{-dual}}}{
\Longleftrightarrow}& \fbox{\eqref{seq:prejet3}}
\end{eqnarray}\end{thm}
\pf Let us denote $\Delta=\Pi^{*}_U$, $\Xi=(\Pi^*_U)^*_{\KStar}$. By
the assumption, we have a $U$-valued pairing
$\conpairing{\cdot,\cdot}_{U}$ of $\Delta$ and $\Pi$ satisfying the
two equalities in (\ref{Eqt:PairingConjugate1}). We also have a
$\KStar$-valued pairing $\conpairing{\cdot,\cdot}_{\KStar}$ of $\Xi$
and $\Delta$ satisfying
\begin{eqnarray*}
\conpairing{\eta,e(\kappa)}_{\KStar}&=& q(\eta)\contract
\kappa,\quad\forall~ \eta\in \Xi,\kappa\in U\otimes\KStar;
\\
\conpairing{r(\zeta),\epsilon}_{\KStar}&=& p(\epsilon)\contract
\zeta,\quad\forall~ \epsilon\in \Delta,\zeta\in V\otimes \KStar.
\end{eqnarray*}
 It suffices to prove that the DVB$^*$ sequences
\eqref{seq:prejet3} and \eqref{seq:prejettrans} are in $V$-duality
and thus it amounts to establish a $V$-valued pairing:
$$\conpairing{\cdot,\cdot}_V:\quad \Xi\times_M \Pi\lon V,$$
such that
\begin{eqnarray}\label{Eqn:tp3}
\conpairing{\eta,i^t(\theta)}_{V}&=& q(\eta)\contract
\theta,\quad\forall~ \eta\in \Xi,\theta\in V\otimes U,
\\\label{Eqn:tp4}
\conpairing{r(\zeta),\sigma}_{V}&=& j(\sigma)\contract
\zeta,\quad\forall~ \epsilon\in \Delta,\zeta\in V\otimes \KStar.
\end{eqnarray}
In fact, the pairing $\conpairing{\eta,\sigma}_V$, for  $\eta\in \Xi,\ \sigma\in\Pi$, can be defined by
\begin{equation}\label{Eqt:conpairingV}
\pairing{\conpairing{\eta,\sigma}_V,v^*} \defbe
\pairing{\conpairing{\eta,\epsilon}_{\KStar},j(\sigma)}
-\pairing{\conpairing{\epsilon,\sigma}_U,q(\eta)},\quad \forall~ \
v^*\in \VStar,
\end{equation}
where an arbitrary $\epsilon\in\Delta$ satisfying
$p(\epsilon)=v^*$ is chosen. This is well defined because if one chooses another $\epsilon'=\epsilon+e(\kappa)$, for some $ \kappa\in
U\otimes \KStar$, then
\begin{eqnarray*}
&&\pairing{\conpairing{\eta,\epsilon'}_{\KStar},j(\sigma)}
-\pairing{\conpairing{\epsilon',\sigma}_U,q(\eta)}\\
&=&\pairing{\conpairing{\eta,\epsilon}_{\KStar},j(\sigma)}
-\pairing{\conpairing{\epsilon,\sigma}_U,q(\eta)}+
\pairing{\conpairing{\eta,e(\kappa)}_{\KStar},j(\sigma)}-\pairing{\conpairing{e(\kappa),\sigma}_U,q(\eta)}\\
&=&\pairing{\conpairing{\eta,\epsilon}_{\KStar},j(\sigma)}
-\pairing{\conpairing{\epsilon,\sigma}_U,q(\eta)}+
\pairing{q(\eta)\contract \kappa,j(\sigma)}-\pairing{j(\sigma)\contract \kappa,q(\eta)}\\
&=&\pairing{\conpairing{\eta,\epsilon}_{\KStar},j(\sigma)}
-\pairing{\conpairing{\epsilon,\sigma}_U,q(\eta)}.
\end{eqnarray*}
We finally prove (\ref{Eqn:tp3}) and (\ref{Eqn:tp4}):
\begin{eqnarray*}
\pairing{\conpairing{\eta,i^t(\theta)}_V,v^*} &=&
\pairing{\conpairing{\eta,\epsilon}_{\KStar},0}
-\pairing{\conpairing{\epsilon,i^t(\theta)}_U,q(\eta)}\\
&=&\pairing{v^*\contract \theta, q(\eta)}=\pairing{q(\eta)\contract
\theta, v^*}.\\ \pairing{\conpairing{r(\zeta),\sigma}_V,v^*} &=&
\pairing{\conpairing{r(\zeta),\epsilon}_{\KStar},j(\sigma)}
-\pairing{\conpairing{\epsilon,\sigma}_U,0}\\
&=&\pairing{v^*\contract \zeta,
j(\sigma)}=\pairing{j(\sigma)\contract \zeta, v^*}.
\end{eqnarray*}
This completes the proof. \qed

\subsection{Duality of associated DVB$^*$ sequences of double vector bundles}

This part of the paper is devoted to show the relevance of duals of
double vector bundles and that of the associated DVB$^*$
 sequences. Given a double vector bundle $(D;A,B;M)_C$ as in Diagram
{\rm(\ref{Fig:DoubleVBDABC})}, we have its dual over $A$, which is
again a double vector bundle, shown as in Diagram
\eqref{Fig:DDualA}. In the meantime, the associated DVB$^*$ sequence
  \eqref{SprayDABC} of $(D;A,B;M)_C$ has an $A^*$-dual, which is denoted by\footnote{This sequence appeared in \cite{alfonso}, where  linear
sections of a double vector bundle are considered.}:
\begin{equation}\label{seq:temp} 0\longrightarrow A^*\otimes
C\longrightarrow \Pi^*_{A^*}\longrightarrow B\longrightarrow 0.
\end{equation}
     We will show that the associated DVB$^*$ sequence of
the dual double vector bundle \eqref{Fig:DDualA} is exactly the
DVB$^*$ sequence \eqref{seq:temp}. For this aim, we first need
another description of the vector bundle $\spray{D}$.

\begin{pro}\label{Pro:SprayExpress} Let
$(\DDualA;A,\CStar;M)_{\BStar}$ be the dual double vector bundle of
$D$ over $A$ given by {\rm(\ref{Fig:DDualA})}. Then,
\begin{equation}\label{Set:DoverA}
(\spray{D})_m \cong \set{\mbox{linear map
}\underline{\sigma}:A_m\longrightarrow (\qDDualACStar)\inverse
(\chi) \mbox{, for some } \chi\in \CStar_m, \mbox{ s.t., }
\qDDualAA\circ \underline{\sigma}=\Id_{A_m}}.
\end{equation}
Under this identification, the maps $i$ and $j$ in
{\rm(\ref{SprayDABC})} are given by
\begin{eqnarray*}
 j&:& \underline{\sigma} \mapsto \chi,\\
i&:& \theta \mapsto (a\mapsto \zeroDDualAA_a \CStarplus
\overline{a\contract\theta},\quad\forall~ a\in A_m),\quad\forall~
\theta\in\AStar_m\otimes\BStar_m\,.
\end{eqnarray*}
\end{pro}
\pf Recall that an element $\sigma\in (\spray{D})_m$ is a
double-linear function on $D_m$.  Notice that $\sigma$
determines a map $\underline{\sigma}$ from $A_m$ to $(\DDualA)_m$ by
sending each $a\in A_m$ to $\underline{\sigma}(a)\in
(\qDDualAA)\inverse(a) $, where
\begin{equation}\label{temp2}
\pairing{\underline{\sigma}(a),d}\defbe \sigma(d),\quad\forall~
(d;a,b;m)\in D_m\,.
\end{equation}
The fact that $\underline{\sigma}(a)\in \DDualA$ is due
$\sigma$ being double-linear. Moreover, we have
\begin{eqnarray*}
\pairing{\underline{\sigma}(a),\zeroDA_a\Bplus \overline{c}}
&=&{\sigma}(\zeroDA_a\Bplus
\overline{c})=\pairing{j(\sigma),c},\quad\forall~ c\in C_m\,.
\end{eqnarray*}
This implies that $\qDDualACStar(\underline{\sigma}(a))=j(\sigma)$
(see Eq.(\ref{Eqt:qDDualACStar})). Hence we obtain a map
$$
\underline{\sigma}: \quad A_m\lon
(\qDDualACStar)\inverse(j(\sigma)),
$$
which is $\Real$-linear and satisfies
$\qDDualAA\circ\underline{\sigma}=\Id_{A_m}$.

Conversely, the element $\underline{\sigma}$  in  (\ref{Set:DoverA})
 determines a double linear function $\sigma$ via relation
(\ref{temp2}). In fact, for any $(d_1;a,b_1;m)\in D_m$, we have
\begin{eqnarray*}
  \sigma(d\Aplus d_1)&=&\pairing{\underline{\sigma}(a),d\Aplus
  d_1}\\
  &=&\pairing{\underline{\sigma}(a),d}+\pairing{\underline{\sigma}(a),d_1}\\
  &=&\sigma(d)+\sigma(d_1),
\end{eqnarray*}
which implies that $\sigma$ is linear with respect to $\Aplus$.
For every $(d_2;a_2,b;m)$, by \eqref{eqt:Cstarplus}, we
have
\begin{eqnarray*}
  \sigma(d\Bplus d_2)&=&\pairing{\underline{\sigma}(a+a_2),d\Bplus
  d_2}\\
  &=&\pairing{\underline{\sigma}(a)\CStarplus \underline{\sigma}(a_2),d\Bplus
  d_2}\\
  &=&\pairing{\underline{\sigma}(a),d}+\pairing{\underline{\sigma}(a_2),d_2}\\
  &=&\sigma(d)+\sigma(d_2),
\end{eqnarray*}
which implies that $\sigma$ is also linear with respect to
$\CStarplus$. Therefore, $\sigma$ is double linear. \qed

Some direct corollaries of Proposition \ref{Pro:SprayExpress} are as follows.
\begin{cor}\label{cor:xd}
Let $(\DDualB;\CStar,B;M)_{\AStar}$ be  the dual double vector
bundle of $D$ over $B$ shown as in Diagram {\rm(\ref{Fig:DDualB})}, then we
have
\begin{equation}\label{Set:DoverA2}
(\spray{D})_m \cong \set{\mbox{linear map
}{\underline{\underline{\sigma}}}:B_m\longrightarrow(\qDDualBCStar)\inverse
(\chi) \mbox{, for some } \chi\in \CStar_m, \mbox{ s.t., }
\qDDualBB\circ \underline{\underline{\sigma}}=\Id_{B_m}}.
\end{equation}
\end{cor}

\begin{cor}\label{cor:sprayDDualAACStar}
Let
\begin{equation}\label{SprayDDualA}
\xymatrix{0 \ar[r] & \AStar\otimes C  \ar[rr]^{i_1} &&
                \spray{\DDualA} \ar[rr]^{j_1} && B \ar[r]  & 0.
                }
\end{equation}
be the associated DVB$^*$ sequence of
$(\DDualA;A,\CStar;M)_{\BStar}$. Then,
$$
(\spray{\DDualA})_m \cong \set{ \mbox{linear map
}\underline{\epsilon}: A_m \longrightarrow(\qDB)\inverse (b) \mbox{,
for some } b\in B_m, \mbox{ s.t., } \qDA\circ
\underline{\epsilon}=\Id_{A_m} }.
$$
Under this identification, the maps $i_1$ and $j_1$ in Sequence
{\rm(\ref{SprayDDualA})} are given by
\begin{eqnarray*}
j_1&:& \underline{\epsilon}\mapsto b,\\
 i_1&:& \kappa\mapsto
(a\mapsto \zeroDA_a\Bplus \overline{a\contract\kappa},\quad\forall
a\in A_m),\quad\forall~ \kappa\in \AStar_m\otimes C_m.
\end{eqnarray*}
\end{cor}
Now we are able to show that the associated DVB$^*$ sequence
\eqref{SprayDDualA} of the double vector bundle $\DDualA$ is
isomorphic to \eqref{seq:temp}, the $A^*$-dual of the associated
DVB$^*$ sequence \eqref{SprayDABC} of $(D;A,B;M)_C$, as the
following theorem claims.
\begin{thm}\label{Thm:ADualCombine}
The associated DVB$^*$ sequence \eqref{SprayDABC} of the double
vector bundle $(D;A,B;M)_C$ and the associated DVB$^*$ sequence
\eqref{SprayDDualA} of the dual double vector bundle
$(\DDualA;A,\CStar;M)_{\BStar}$ are in $\AStar$-duality.
\end{thm}
\pf  We need to establish an $\AStar$-pairing
$\conpairing{\cdot,\cdot}_{\AStar}$ between $\spray{D}$ and
$\spray{\DDualA}$ satisfying the following two equalities:
\begin{eqnarray*}
 \conpairing{i(\theta),\epsilon}_{\AStar}&=&j_1(\epsilon\contract\theta),
 \quad \forall~\epsilon\in \spray{\DDualA},\quad \theta\in \AStar_m\otimes
 \BStar_m,\\
\conpairing{\sigma,i_1(\kappa)}_{\AStar}
&=&j(\sigma)\contract\kappa,\quad \forall~\sigma\in \spray{D},\quad
\kappa\in \AStar_m\otimes C_m.
\end{eqnarray*}
This can be done by defining $\conpairing{\sigma,\epsilon}_{\AStar}$,
for $\sigma\in \spray{D}_m$ and $\epsilon\in
\spray{\DDualA}_m$ to be
$$
\conpairing{\sigma,\epsilon}_{\AStar}(a)\defbe
\pairing{\underline{\sigma}(a),\underline{\epsilon}(a)},
\quad\forall~ a\in A_m\,.
$$
It is obviously bilinear and symmetric. Moreover, Proposition
\ref{Pro:SprayExpress} and Corollary  \ref{cor:sprayDDualAACStar} imply that $$
\conpairing{\sigma,\epsilon}_{\AStar}(a)=\sigma(\underline{\epsilon}(a))=\epsilon(\underline{\sigma}(a)).
$$
Thus, by definition of $i(\theta)$ in
(\ref{Eqt:AStarBStarFunctionOnSlice}), we have
\begin{eqnarray*}
\conpairing{i(\theta),\epsilon}_{\AStar}(a) &=&
{i(\theta)}(\underline{\epsilon}(a))= \pairing{\theta,a\otimes
j_1(\epsilon)}=(j_1(\epsilon)\contract \theta) (a).
\end{eqnarray*}
Similarly, the map $i_1$ in (\ref{SprayDDualA}) is given by
$$
i_1(\kappa)(\Phi)=\pairing{\kappa,\qDDualAA(\Phi)\otimes\qDDualACStar(\Phi)},
\quad\forall~ \Phi\in \DDualA.
$$
So we obtain
\begin{eqnarray*}
\conpairing{\sigma,i_1(\kappa)}_{\AStar}(a) &=&
i_1(\kappa)(\underline{\sigma}(a))= \pairing{\kappa,a\otimes
j(\sigma)}=(j(\sigma)\contract\kappa)(a).
\end{eqnarray*}
This completes the proof. \qed

Since the categories  {\bf DVB}  and {\bf DVBS} are equivalent (Theorem \ref{Thm:DVB1-1CM}), a direct corollary is the following:
\begin{cor}\label{cor:DEpair}
Suppose that $(D;A,B;M)_C$ and $(E;A,\CStar;M)_{\BStar}$ are double
vector bundles with a side bundle $A$ in common and with cores $C$
and $\BStar$ respectively. Suppose that their associated DVB$^*$
sequences
$$
\xymatrix@C=0.5cm{
  0 \ar[r] & \AStar\otimes \BStar \ar[rr]^{i} &&
\spray{D}
   \ar[rr]^{j} && \CStar \ar[r] & 0 },
$$
$$
\xymatrix@C=0.5cm{
  0 \ar[r] & \AStar\otimes C \ar[rr]^{i_1} &&
\spray{E}
   \ar[rr]^{j_1} && B \ar[r] & 0 }.
$$
are in $\Astar$-duality as  DVB$^*$ sequences. Then there induces a
paring between $D$ and $E$ over $A$ so that they are mutually duals
as double vector bundles .
\end{cor}

We also recover a well known fact  \cite[Theorem 9.2.2]{Mkz:GTGA}.
\begin{cor}\label{cor:dualoverCstar}
The double vector bundles $\DDualA$ and
$\DDualB$ are mutually duals over $C^*$.
\end{cor}
\pf By Theorem \ref{Thm:TransLDual} and \ref{Thm:ADualCombine}, the
associated DVB$^*$ sequences of $\DDualA$ and $\DDualB$ are in
$C$-duality. Then the conclusion follows directly from Corollary
\ref{cor:DEpair}. \qed \vspace{3mm}

\section{The Atiyah algebroid and the jet bundle.}
Let $E\stackrel{q}{\longrightarrow}M$ be a vector bundle and
$\EStar\stackrel{\qStar}{\longrightarrow}M$  the dual bundle
of $E$. The tangent double vector bundle of $E$ is a double vector
bundle $(\TE;TM,E;M)_{E}$ which fits the following diagram:
\begin{equation}\label{Fig:TE}
\xymatrix{
  \TE \ar[d]_{\Tq} \ar[r]^{\pE}
                & E \ar[d]^{q} &   \\
  TM  \ar[r]_{p}
                & M & \ar[l]^{q} E.        }
\end{equation}
Here $TE\stackrel{\Tq}{\longrightarrow}{TM}$ is the   tangent to the   $E\stackrel{q}{\longrightarrow}M$ and
$TE\stackrel{\pE}{\longrightarrow}{E}$ is the usual tangent
bundle. We denote elements of $\TE$ by $\xi,\eta,\zeta$
... and we write $(\xi;x,e;m)$ to indicate that $e=\pE(\xi)$,
$x=\Tq(\xi)$, and $m=p\circ\Tq(\xi)=q\circ\pE(\xi)$. With respect to
the tangent bundle structure $(TE,\pE,E)$, we use usual notations.
The zero element in $T_eE$ is denoted by $\tilde{0}_e$. In the
prolonged tangent bundle structure, $(TE,\Tq,TM)$, we use $\TMplus$
for addition, $\TMminus$ for subtraction, and $\TMtimes$ for scalar
multiplication. The zero element of the fiber $(\Tq)\inverse(x)$ is
denoted by $T(0)(x)$.
\vskip 0.5cm
The embedding of the core $\overline{(\cdot)}$: $E\lon \TE$ is to
identify the tangent space $T_{0_m}(E_m)$ canonically with $E_m$.
That is, the element $\overline{e}$ of $T_{0_m}(E_m)$ corresponding
to $e\in E_m$ is given by
$$
\overline{e}=\frac{d}{dt}|_{t=0}(te).
$$

The vertical dual of (\ref{Fig:TE}) is actually the tangent double
vector bundle of $\EStar$:
\begin{equation}\label{Fig:TEStar}
\xymatrix{
  \TEStar \ar[d]_{\TqStar} \ar[r]^{\pEStar}
                & \EStar \ar[d]^{\qStar} &   \\
  TM  \ar[r]_{q}
                & M & \ar[l]^{\qStar} \EStar.        }
\end{equation}

Dualizing the   double vector bundle $TE$ in Diagram \eqref{Fig:TE}  over $E$ leads to a double vector bundle
$(\TStarE;E,\EStar;M)_{\TStarM}$:
\begin{equation}\label{Fig:TStarE}
\xymatrix{
  \TStarE \ar[d]_{\cE} \ar[r]^{\rE}
                & \EStar \ar[d]^{\qStar} &   \\
  E  \ar[r]_{q}
                & M & \ar[l]^{\pStar} \TStarM,        }
\end{equation}
which will be referred to as the cotangent double vector bundle of $E$,
or cotangent dual of $\TE$. An element of $\TStarE$ will be denoted by
$(\frke;e,\varphi;m)$ to indicate that $e=\cE(\frke)$,
$\varphi=\rE(\frke)$, $m=q(e)=\qStar(\varphi)$. The map $\rE:
\TStarE\lon \EStar$ is defined as follows. Take $\frke\in
T_{e}^*(E)$, where $e\in E_m$. Define
$$
\pairing{\rE(\frke), e'}=\pairing{\frke,\tilde{0}_e \TMplus
\overline{e'}},\quad \forall~ e'\in E_m\,.
$$
The embedding of $w\in T^*_mM$ into $\TStarE$ is given by
$$
\pairing{\overline{w},T(0)(x)+\overline{e}}=\pairing{w,x},\quad
\forall e\in E_m, x\in T_mM,\quad\mbox{( i.e. } \overline{w}=q^*(w)
\mbox{ ).}
$$

There is of course an analogously
\begin{equation}\label{Fig:TStarEStar}
\xymatrix{
  \TStarEStar \ar[d]_{\cEStar} \ar[r]^{\rEStar}
                & E \ar[d]^{q} &   \\
  \EStar  \ar[r]_{\qStar}
                & M & \ar[l]^{\pStar} \TStarM,        }
\end{equation}
which is the cotangent double vector bundle of $\EStar$.

\vskip0.2in \noindent\textbf{$\bullet$The Atiyah algebroid $\dev E$}\vspace{3mm}

For a vector bundle $E\stackrel{q}{\lon}M$,   we denote the gauge Lie algebroid of the
 frame bundle
 $F(E)$ by
$\dev E$, which is also called the Atiyah algebroid, also known as the covariant differential operator bundle of $E$ (see \cite[Example 3.3.4]{Mkz:GTGA}
 and \cite{MackenzieX:1994}).   Here we treat each
element $\frkd$ of $\dev{E}$ at $m\in M$ as an $\Real$-linear
operator $\Gamma(E)\lon E_m$   with some $x\in T_mM$ (which
is uniquely determined by $\frkd$ and called the anchor of $\frkd$)
such that
$$
\frkd(fu)=f(m)\frkd(u)+x(f)u(m), \, \quad \quad ~~ \forall~
f\in\CWM,
 u\in\Gamma(E).
$$
It is  known that $\dev{E}$ is a  transitive Lie algebroid over $M$
\cite{KSK:2002}. The anchor of $\dev{E}$ is given by $\jd(\frkd)=x$
and the Lie bracket $[\cdot,\cdot ]_{\dev}$ of $\Gamma(\dev{E})$ is
just the  commutator.

As a transitive Lie algebroid, $\dev{E}$ has an associated exact sequence, called the  Atiyah  sequence:
\begin{equation}\label{Seq:DE}
\xymatrix@C=0.5cm{0 \ar[r] & \EStar\otimes E  \ar[rr]^{\id} &&
                \dev{E}  \ar[rr]^{\jd} && TM \ar[r]  & 0.
                }
\end{equation}

We have an equivalent description of
$\dev{E}$:
$$
(\dev{E})_m=\set{(\frkd,x) \mbox{ where } x\in T_mM, \frkd\in
\Hom(E^*_m,(T(q_*))\inverse(x)) \mbox{ and }
p_{\EStar}\circ\frkd=\Id_{E^*_m}}.
$$
From this point of view, a section of $\dev{E}$, namely a derivation
$(\frkd,x)$, is a linear vector field $\frkd$ on $\EStar$ which
projects to a vector field $x\in\frkX(M)$, such that for any $u\in\Gamma(E)$, if it is
treated as a fiber-wise linear function $l_u$, then
$$
\frkd(l_u)=l_{\frkd(u)}.
$$
One may refer to \cite[Definition
3.4.1]{Mkz:GTGA} for more details.

By Corollary
\ref{cor:xd},  we are able to  draw the following conclusion.
\begin{cor}\label{cor:devE}
The exact sequence \eqref{Seq:DE} associated with the Atiyah
algebroid $\dev E$ is the associated DVB$^*$ sequence of the
cotangent double vector bundle of $E^*$ shown as in Diagram
{\rm(\ref{Fig:TStarEStar})}, i.e.
$$
\dev{E}=\spray{\TStarEStar}.
$$
\end{cor}

\vskip0.2in \noindent\textbf{$\bullet$ The Jet bundle $\jet
E$}\vspace{3mm}

Given a  vector bundle
$E\stackrel{q}{\lon}M$ and a pint $m\in M$,
for any $e\in E_m$, the tangent space $T_e(E_m)$, called the vertical subspace, is a linear
subspace of the full tangent space of $E$ at $e$.
  The full tangent space can be
decomposed into a direct sum of $T_e(E_m)$ and a
complementary horizontal subspace. We can define a fiber bundle
$\jet E$ over $E$ whose fiber at $e$ is the set of all possible
horizontal subspaces:
$$
(\jet{E})_e\defbe\set{\mbox{linear map } \mu: T_mM\lon T_eE,\mbox{
such that } \Tq\circ\mu=\Id_{T_mM}}.
$$
There is a linear structure of $\jet E$  as follows: for  $\mu_1\in (\jet{E})_{e_1}$, $\mu_2\in (\jet{E})_{e_2}$ with
$q(e_1)=q(e_2)=m$, $r\mu_1+\mu_2\in (\jet E)_{r e_1+e_2}$ (where $r\in\mathbb{R}$) is a linear map $
 T_mM\lon T_{re_1+e_2}E $ defined by:
$$
(r\mu_1+\mu_2)(x)= r\TMtimes\mu_1(x)\TMplus\mu_2(x),\quad\forall~
x\in T_mM.
$$

Equipped with this operation, $\jet E$ is a vector
bundle over $M$, called the \emph{first order jet bundle} over $M$,
or simply the jet bundle of $E$ (see \cite{Saunders} for more details about jet bundles) . In summary, the fiber of $\jet{E}$
at $m$ is the collection
$$
\jet{E}_m=\set{(\mu,e),\ \mbox{where } e\in E_m, \mu\in \Hom(T_mM,
T_eE)\mbox{ such that } \Tq\circ\mu=\Id_{T_mM}},
$$
with the linear structure defined by
$$
r(\mu_1,e_1)+(\mu_2,e_2)=(r\mu_1+\mu_2,re_1+e_2),
$$
for $(\mu_p,e_p)\in\jet{E}_m$ ($p=1,2$), $r\in\Real$. It is well
known that there is an associated exact sequence:
\begin{equation}\label{Seq:JetE}
\xymatrix@C=0.5cm{0 \ar[r] & \TStarM\otimes E  \ar[rr]^{\e} &&
                \jet{E} \ar[rr]^{\p} && E \ar[r]  & 0.
                }
\end{equation}

By
Proposition \ref{Pro:SprayExpress}, we have

\begin{cor}\label{cor:SprayJetOfE}
The exact sequence \eqref{Seq:JetE} associated with $\jet E$ is the
associated DVB$^*$ sequence of the tangent double vector bundle of
$\EStar$ shown as in Diagram {\rm(\ref{Fig:TEStar})}, i.e.
$$\jet{E}=\spray{\TEStar}.$$
\end{cor}

\vskip0.2in \noindent\textbf{$\bullet$ The relations between $\dev E$, $\jet
E$, $\dev E^*$ and $\jet E^*$}\vspace{3mm}

Since $TE$ and $TE^*$ are double vector bundles in duality over
$TM$, and $\jet{E^*}=\spray{\TE}$, $\jet{E}=\spray{\TEStar}$
(Corollary \ref{cor:SprayJetOfE}),
 we have an induced $T^*M$-pairing between $\jet {E^*}$ and $\jet E$ (Theorem \ref{Thm:ADualCombine}). In other words, they are mutually $T^*M$-duals as  DVB$^*$ sequences.

Similarly, $T^*E^*$ and $TE^*$ are in duality over $E^*$ as double vector bundles, and $\jet{E}=\spray{\TEStar}$, $\dev E=\spray{T^*E^*}$(Corollary \ref{cor:devE}),
 we come to the conclusion that $\jet{E}$ and $\dev E$ are mutually $E$-duals as  DVB$^*$ sequences. The pairing between $\jet E$ and $\dev E$ is exactly the
 one defined in \cite{CLomni}.  For the same reasons, $\jet{E^*}$ and $\dev E^*$ are mutually $E^*$-duals.

In \cite[Lemma 2.3]{CLomni}, it is shown that $\dev E$ is
canonically isomorphic to $\dev E^*$. The isomorphism is just the
transposition of a  DVB$^*$ sequence. Thus, by Theorem
\ref{Thm:DVB1-1CM} and Corollary \ref{cor:devE}, we recover a well
known result:
\begin{cor}
The double vector bundle associated to $T^*E$, is isomorphic to the double vector bundle
associated to $T^*E^*$.
\end{cor}

Recall Diagram \eqref{Dig:square}, we can again  summarize the relations between $\dev E$, $\jet
E$, $\dev E^*$ and $\jet E^*$ by the following diagram:
\begin{eqnarray}\nonumber \dev E^* &\quad\stackrel{\mbox{\emph{ transposition}}}{
\Longleftrightarrow}\quad& \dev E\\\nonumber
\\\nonumber E^* \mbox{\emph{-dual}}{\bigg \Updownarrow} &
& {\bigg \Updownarrow}E
\mbox{\emph{-dual}}\\\nonumber\quad\\\nonumber
\jet E^* &\stackrel{T^*M\mbox{\emph{-dual}}}{
\Longleftrightarrow}& \jet E.
\end{eqnarray}

\end{document}